\documentclass[letterpaper,11pt]{article}
\usepackage[letterpaper, left=1in,right=1in,top=1in,bottom=1in]{geometry}   
\usepackage{color}
\usepackage{amssymb}
\usepackage{epsfig}
\usepackage{amsmath}
\usepackage{amsfonts}
\usepackage{graphicx}
\usepackage{subfigure}
\usepackage{color}
\usepackage{graphicx}
\usepackage{epsfig}
\usepackage{cite}
\usepackage{amsmath}
\usepackage{amssymb}
\usepackage{mathrsfs}
\usepackage{amsfonts}
\usepackage{dsfont}
\usepackage{graphicx}
\usepackage{mathrsfs}
\usepackage{hyperref}
\usepackage{hyperref}
\usepackage{titlesec}
\usepackage[T1]{fontenc}
\usepackage{ragged2e}
\justifying

\DeclareGraphicsExtensions{.pdf,.png,.jpg,.JPG}

\usepackage{lastpage}
\usepackage{fancyhdr}
\raggedright
\titleformat{\section}{
	\vspace{12pt}\scshape\raggedright\large
}{}{0em}{\textsc}

\begin{document}
	\thispagestyle{empty}
	
	\begin{tabular*}{\textwidth}{l@{\extracolsep{\fill}}r}
		\textbf{\centerline{A Note for CPS Data-driven Approaches Developed in the IDS Lab}}\vspace{5pt}\\ 
		\centerline{Andreas A. Malikopoulos}\\
        \centerline{Cornell University}\\
        \centerline{amaliko@cornell.edu}\\
	\end{tabular*}

\vspace{+5EX}

\centerline{\textbf{Abstract}}

\justify{

The rapid evolution of Cyber-Physical Systems (CPS) across various domains like mobility systems, networked control systems, sustainable manufacturing, smart power grids, and the Internet of Things  necessitates innovative solutions that merge control and learning \cite{Malikopoulos2022a}. Traditional model-based control methodologies often fail to adapt to the dynamism and complexity of modern CPS. This report outlines a comprehensive approach undertaken by the Information and Decision Science (IDS) Lab, focusing on integrating data-driven techniques with control strategies to enhance CPS performance, particularly in the context of energy efficiency and environmental impact.
CPS are intricate networks where physical and software components are deeply intertwined, operating as systems of systems. These systems are characterized by their informationally decentralized nature, posing significant challenges in optimization and control. Classical control methods depend heavily on precise models, which often do not capture the full complexity of real-world CPS. As these systems generate large volumes of real-time data, there is a growing need for control algorithms that can leverage this data effectively. The IDS Lab is at the forefront of developing such data-driven approaches for CPS.
	
	\section{\textbf{Autonomous Powertrain Systems}}
	\justify{
		
		We have conducted research on autonomous powertrain control systems capable of learning and adapting to drivers’ behaviors to optimize engine performance and fuel efficiency \cite{Malikopoulos2011b, Malikopoulos2007,Malikopoulos2007a,Malikopoulos2009b, Malikopoulos2009a,Malikopoulos2006,Malikopoulos2011a,Malikopoulos2015,Malikopoulos2015b,Malikopoulos2008b,Malikopoulos2010a,Malikopoulos2010,Malikopoulos2009e,Malikopoulos2010b, Malikopoulos2009,Malikopoulos2013d}. This interest was sparked by an article highlighting the significant impact of driving styles on vehicle performance. We have developed a theoretical framework for transforming a vehicle’s engine into an intelligent system that learns its optimal operation dynamically. By modeling the driver’s behavior as a control Markov chain, the engine learns and adapts to achieve optimal performance criteria, such as minimizing fuel consumption and emissions. This framework’s practical applications led to a US patent and subsequent developments at General Motors, where I contributed to creating self-learning control algorithms for advanced powertrain systems.

		We have also broadened our research to include stochastic optimal control and power management for hybrid-electric vehicles  and plug-in hybrid-electric vehicles  	\cite{Park2011,Malikopoulos2011,Rios-Torres2016a,Malikopoulos2013a,Malikopoulos2016,Malikopoulos2013b,Malikopoulos2012,Malikopoulos2013,Malikopoulos2014,Malikopoulos2015a,Malikopoulos2014a,Shaltout2014,Pourazarm2014, Shaltout2015b,Malikopoulos2015d,Malikopoulos2016b}. Our work involved developing algorithms for optimal energy distribution and power management in these vehicles, enhancing their efficiency and sustainability. In addition, we have investigated the integration of these vehicles with smart grids and the implications of their widespread adoption on urban infrastructure \cite{Sharma2016,Dong2017a,Dong2016}. My role as Deputy Director of the Urban Dynamics Institute at Oak Ridge National Laboratory allowed me to spearhead initiatives that examined the environmental impacts of conected and automated vehicles (CAVs) and sought to improve urban transportation systems’ sustainability and accessibility.
	}

\section{\textbf{Emerging Mobility Systems}}
\justify{
	
Emerging mobility systems \cite{chremos2020MobilityMarket,chremos2021MobilityGame,chremos2020SharedMobility} represent a new paradigm in transportation, integrating advanced technologies and control algorithms to enhance efficiency and sustainability 	\cite{chalaki2021Reseq, Ray2021DigitalCity,mahbub2021_platoonMixed,Sumanth2021,mahbub2022_ifac, Malikopoulos2019CDC, Mahbub2019ACC,  malikopoulos2019ACC, Malikopoulos2020, mahbub2020sae-1, mahbub2020ACC-2, chalaki2020hysteretic, Zhao2018CTA, tzortzoglou2023approach, chalaki2021CSM,mahbub2020Automatica-2, jang2019simulation, mahbub2020decentralized, Zhao2019CCTA-1, mahbub2023_automatica, Connor2020ImpactConnectivity, mahbub2022ACC,Le2022CDC, zhao2019enhanced, tzortzoglou2024feasibility, Le2023CDC, Malikopoulos2018,K122021AExperiment,Malikopoulos2017,Le2023Stochastic, Malikopoulos2018b,mahbub2020sae-2, le2024controller, Rios-Torres2015,Rios-Torres2017b,Rios-Torres2017a,Rios-Torres2017,Rios-Torres2016b,Zhao2018ITSC, tzortzoglou2023performance,zhao2019CCTA-2,  Rios-Torres2018,Stager,Zhao2018,Zhang2017a,Zhang,Zhao2018a, assanis2018ITSC,Beaver2020DemonstrationCity,chalaki2019optimal,chalaki2020experimental,chalaki2020ICCA,Mahbub2020ACC-1, chalaki2020TCST,Nishanth2023AISmerging,chalaki2020TITS}. These systems include CAVs, shared mobility platforms, and interconnected public transit networks. The IDS Lab's research in this area has focused on developing control algorithms for the optimal coordination of CAVs, addressing challenges such as collision avoidance, traffic flow optimization, and energy efficiency. One significant advancement was the development of a decentralized optimal control framework for CAVs at merging roadways, which significantly improves fuel efficiency and reduces emissions by optimizing vehicle trajectories in real-time.	
	
Autonomous Mobility-on-Demand (AMoD) systems \cite{bang2024confidence, Bang2023flowbased,bang2021AEMoD, bang2023exploring,bang2024cts} represent a significant evolution in urban transportation \cite{Chremos2020MechanismDesign}, offering dynamic, real-time routing \cite{Bang2022rerouting, Bang2022combined,chremos2023AtomicRouting,bang2023optimal} and scheduling of autonomous vehicles. The IDS Lab's research in this area has focused on developing algorithms for the efficient operation of electric AMoD systems, including routing optimization and charging infrastructure planning. By integrating CAV behavior into AMoD systems, we have developed hierarchical optimization frameworks that improve traffic flow and system scalability. This research addresses the technical challenges of AMoD systems while also considers the broader impacts on urban mobility, sustainability, and equity \cite{bang2024emergingequity,bang2024cts,Bang2023mem}.
}

\section{\textbf{Human Involvement in Cyber-Physical Systems}}
\justify{
Considering human factors in the design and operation of CPS is crucial, particularly in socio-technical systems where human behavior significantly impacts system performance. The IDS Lab's in this area research has explored the interactions between AI-driven systems and human users, developing algorithms that adapt to human behaviors to optimize system performance \cite{dave2024airecommend,faros2023adherence}. Additionally, I have proposed novel metrics for assessing mobility equity in transportation networks, ensuring that advancements in CPS contribute to social equity and inclusivity.

}

\section{\textbf{Team Theory}}
\justify{

Team theory offers a valuable framework for addressing decentralized control problems involving multiple agents with shared objectives but different information. The IDS Lab's research in this domain has focused on sequential dynamic team decision problems, providing structural results and a dynamic programming decomposition for such problems \cite{Malikopoulos2024,Malikopoulos2022a,Malikopoulos2021,Dave2021minimax,dave2019decentralized,Dave2020,Dave2021nestedaccess,Dave2020a,Dave2021a}. These contributions have important implications for the design and control of informationally decentralized systems like CPS, enabling more efficient and robust coordination among system components \cite{Nishanth2023AISmerging,Dave2023infhorizon,venkatesh2023stochastic,dave2022additive,Dave2023approximate, Dave2022approx}.

}

\section{\textbf{Mechanism Design Theory}}
\justify{

Mechanism design theory, a subfield of game theory, provides a structured approach to designing rules and protocols that lead to desired outcomes, even when participants have private information and individual incentives. The IDS Lab's research in this area has focused on applying mechanism design to CPS, particularly in contexts where agents need to make decisions that are socially optimal \cite{chremos2022CSMArticle,Dave2020SocialMedia}.
One major application of mechanism design in my research has been in the development of incentive schemes for energy-efficient behavior in mobility markets \cite{chremos2023AtomicRouting}. By designing mechanisms that align individual incentives with collective goals, we can encourage behaviors that enhance system-wide efficiency and sustainability.

}

\section{\textbf{Multi-Agent System}}
\justify{

Multi-agent systems  play a critical role in the landscape of CPS, where numerous autonomous entities interact and collaborate to achieve complex tasks \cite{Beaver2020AnFlockingb}. The IDS Lab's research in this area has focused on developing robust coordination algorithms that ensure seamless interaction among agents, which may include vehicles, drones, or other autonomous systems  \cite{Beaver2020AnAgents,Beaver2020Energy-OptimalConstraints,bang2021energy,Beaver2019AGeneration,Beaver2021Constraint-DrivenStudy,Beaver2020BeyondFlocking,Beaver2020AnFlocking}. These algorithms address critical challenges such as task allocation, resource sharing, conflict resolution, and cooperative decision-making. One notable project involved designing a decentralized control framework for multi-agent systems, enabling agents to operate with limited information about the global state while still achieving collective goals. This approach leverages principles of game theory and distributed optimization to ensure that agents can make informed decisions based on local observations and interactions. In the context of emerging mobility systems, MAS research has been pivotal in optimizing the behavior of connected and automated vehicles within mixed traffic environments. By simulating various scenarios, we have developed strategies that enhance traffic flow, reduce congestion, and improve safety. These strategies involve real-time data exchange among vehicles and infrastructure, enabling adaptive responses to dynamic traffic conditions. 
Furthermore, my research has explored the integration of MAS with machine learning techniques, particularly reinforcement learning, to enable agents to learn and adapt to their environment over time. This combination of MAS and learning algorithms has shown promise in applications ranging from urban air mobility to automated logistics and smart manufacturing.

}

\newpage

\bibliographystyle{IEEEtran}
\bibliography{references1,IDS_Publications_06132024,TAC_Ref_structure,TAC_Ref_Andreas}

\end{document}